\documentclass[journal]{IEEEtran}
\usepackage{cite}
\usepackage{amsmath}
\renewcommand{\vec}[1]{\boldsymbol{#1}}
\usepackage{setspace}
\usepackage{bm}
\usepackage{amssymb}
\usepackage{latexsym}
\usepackage{makecell}
\usepackage{graphicx}
\usepackage{url}
\usepackage{tabularx}
\usepackage{multirow}
\usepackage{float}
\usepackage{slashbox}
\usepackage{algorithm}
\usepackage{algorithmic}
\usepackage[subfigure]{graphfig}

\usepackage{color}

\begin{document}

\title{Distributed Optimal Gas-Power Flow Using \\ Convex Optimization and ADMM}
\author{Cheng~Wang,~\IEEEmembership{Student Member,~IEEE}, Wei~Wei,~\IEEEmembership{Member,~IEEE}, Jianhui~Wang,~\IEEEmembership{Senior Member,~IEEE}, Linquan~Bai,~\IEEEmembership{Student Member,~IEEE}, Yile~Liang,~\IEEEmembership{Student Member,~IEEE}
\thanks{The work of W. Wei is supported by the National Natural Science Foundation of China (51321005, 51577163). The work of J. Wang is sponsored by the U. S. Department of Energy. (\emph{Correspond author: Wei Wei})} 
\thanks{C. Wang, W. Wei, and Y. Liang are with the State Key Laboratory of Power Systems, Department of Electrical Engineering, Tsinghua University, 100084 Beijing, China. (e-mail: c-w12@mails.tsinghua.edu.cn, wei-wei04@mails.tsinghua.edu.cn, yileat17@163.com).}
\thanks{J. Wang is with Argonne National Laboratory, Argonne, IL 60439, USA (e-mail: jianhui.wang@anl.gov).}
\thanks{L. Bai is with Department of Electrical Engineering and Computer Science, University of Tennessee, Knoxville, USA (e-mail: lbai3@vols.utk.edu)}}
\maketitle

\begin{abstract}
— This paper proposes a convex optimization based distributed algorithm to solve multi-period optimal gas-power flow (OGPF) in coupled energy distribution systems. At the gas distribution system side, the non-convex Weymouth gas flow equations is convexified as quadratic constraints. The optimal gas flow (OGF) subproblem is solved by an iterative second-order cone programming procedure, whose efficiency is two orders of magnitudes higher than traditional nonlinear methods. A convex quadratic program based initiation scheme is suggested, which helps to find a high-quality starting point. At the power distribution system side, convex relaxation is performed on the non-convex branch flow equations, and the optimal power flow (OPF) subproblem gives rise to a second order cone program. Tightness is guaranteed by the radial topology. In the proposed distributed algorithm, OGF and OPF are solved independently, and coordinated by the alternating direction multiplier method (ADMM). Numerical results corroborate significant enhancements on computational robustness and  efficiency compared with existing centralized  OGPF methods.
\end{abstract}


\section*{Nomenclature}
Most of the symbols and notations used throughout this paper are defined below for quick reference. Others are defined following their first appearances, as needed.

\subsection{Sets and Indices}
\begin{IEEEdescription}[\IEEEusemathlabelsep\IEEEsetlabelwidth{$aaaaaaaa$}]
\item [$t\in T$]  Time periods
\item [$g\in G$]  Non-gas distributed generator (DG)
\item [$n\in N$]  Gas-fired DG
\item [$i_p \in I_p$]   Power distribution network (PDN) nodes
\item [$l_p \in L_p$]   PDN lines
\item [$d_p \in D_p$]   PDN loads
\item [$w \in W$]   Gas retailers
\item [$c \in C$]   Gas compressors (Gas active pipelines)
\item [$i_g \in I_g$]   Gas distribution network (GDN) nodes
\item [$l_g \in L_g$]   Gas passive pipelines
\item [$d_g \in D_g$]   GDN loads
\end{IEEEdescription}

\subsection{Parameters}
\begin{IEEEdescription}[\IEEEusemathlabelsep\IEEEsetlabelwidth{$aaaaaaaaaa$}]
\item [$P_{g}^{min}/P_{g}^{max}$]  Active power range of non-gas DGs
\item [$P_{n}^{min}/P_{n}^{max}$]  Active power range of gas-fired DGs
\item [$Q_{g}^{min}/Q_{g}^{max}$]  Reactive power range of non-gas DGs
\item [$Q_{n}^{min}/Q_{n}^{max}$]  Reactive power range of gas-fired DGs
\item [$I_{l_p}^{max}$]  PDN line capacity
\item [$p_{d_pt}$]  PDN load demand
\item [$r_{l_p}/x_{l_p}$]  Resistance/reactance of PDN line
\item [$Q_{g}(\cdot)$]  Generation cost of non-gas DGs
\item [$\tau^u_{i_g}/\tau^l_{i_g}$]  Gas pressure range
\item [$\beta_n$]   Gas-electricity conversion factor
\item [$\gamma_c$]  Compression factor of the compressor
\item [$\phi_{l_g}$]  Weymouth equation coefficient
\item [$\alpha_c$]  Fuel consumption coefficient
\item [$Q_{wt}$]  Gas Price
\item [$q_{d_gt}$]  Gas distribution network (GDN) load
\item [$T_k$]   Temperature
\item [$Z_{l_g}$]  Compression factor of the pipeline
\item [$x_{l_g}$]  Length of the pipeline
\item [$R_{l_g}$]  Diameter
\item [$\rho_0$]  Gas density in standard condition
\item [$\mu$]  Specific gas constant
\item [$F_{l_g}$]  Pipeline friction coefficient
\item [$\lambda$]  Unit transformation constant
\item [$\chi$]  Thermal equivalent conversion constant
\end{IEEEdescription}

\subsection{Variables}
\begin{IEEEdescription}[\IEEEusemathlabelsep\IEEEsetlabelwidth{$aaaaaaaa$}]
\item [$p_{gt},p_{nt}$]  Active power of DGs
\item [$q_{gt},q_{gt}$]  Reactive power of DGs
\item [$pf_{l_pt}/qf_{l_pt}$] Active/reactive power of PDN lines
\item [$\nu_{i_pt}$] Nodal voltage square of PDN
\item [$I_{i_pt}$] Line current square of PDN
\item [$y_{wt}$]  Purchased gas
\item [$y_{l_gt}^{in}/y_{l_gt}^{out}$]  Gas in/out flow of passive pipeline
\item [$y_{ct}^{in}/y_{ct}^{out}$]  Gas in/out flow of active pipeline
\item [$u_{i_gt}$]  Nodal pressure of GDN
\item [$m_{l_gt}$]  Average gas mass of GDN
\end{IEEEdescription}

\section{Introduction}

\IEEEPARstart{I}n the past decades, thanks to the dramatic decrease in natural gas price benefiting from the shale rock revolution \cite{ShaleGas}, as well as breakthroughs in turbine technology, gas-fired units are playing a more important role in power systems. Besides the economic profits and environmental benefit \cite{Mohammad_Gas}, physical connections between the gas system and power system create remarkable interdependency across these critical infrastructures, from transmission level to distribution level. The fuel adequacy of gas-fired units impacts the security of both networks and has been studied in \cite{Mark_Fuel} and \cite{Chertkov_Fuel}. Joint operation of power and natural gas system are desired. Many interesting and inspiring studies have been found focusing on this coordination issues. A robust co-scheduling model of a gas-electric system is proposed in \cite{Cong_Look_ahead_congestion}, which takes the gas congestion into account. A bi-level interdiction model is formulated in \cite{Khodayar_gas_AD} to detect the most severe attack to the gas-electric system in the distribution level. The long-term independency of natural gas and power system is considered in the co-planning model proposed in \cite{Xiaping_Planning}.

Optimal gas-power flow (OGPF) is the most fundamental problem in the coordinated system operation, and has been widely discussed. For the power system model, most existing work adopts the direct current (DC) power flow model for simplicity \cite{Chaudry_DCOPF,Cong_Dispatch,Cong_PDE,Anatoly_PDE,Chengcheng_MILP}. The accurate alternating current (AC) power flow formulation is used in \cite{Alberto_Newton,Sheng_Probabilistic}, however, non-convex AC power flow model makes the computation much more challenging. For the natural gas system model, the steady-state assumption is widely adopted, i.e., the in-flow and out-flow of each pipeline are equal\cite{Chengcheng_MILP,Alberto_Newton, Wolf_BNG,Carlos_MILP}. The Weymouth equation, which describes the relationship between gas flow and nodal pressure, is also non-convex and imposes a major challenge on seeking the optimal solution. A common treatment is the mixed integer linear program (MILP) based piecewise linear approximation, such as that in \cite{Chengcheng_MILP} and \cite{Carlos_MILP}, where the accuracy can be well-controlled by adjusting the number of linear segments. However, the computation efficiency remains a major challenge as MILPs are NP-hard. Another choice is the Newton's method \cite{Alberto_Newton}, in which a feasible initial point is required. Certainly, nonlinear programming solvers, such as IPOPT, are able to handle nonlinear equations directly \cite{Anatoly_PDE}. Though the existing works about OGPF problem mainly focus on transmission level, the interdependency between power and gas system on distribution level is stronger \cite{GE_GasPower} and getting more attention \cite{Urban_GasPower,Distribution_GasPower}, indicating the necessity of distribution-level OGPF analysis.

In consideration of relatively slow gas dynamics as well as the long distance of gas pipelines, non-steady-state operating condition is usual in natural gas system, and it is reasonable to model system dynamics, referred to as the linepack. Non-steady-state operation endows the gas system   with storage capability which also adds flexibility to power system operation  \cite{Nico_SingleDirection}. Meanwhile, it is observed in \cite{Anatoly_PDE} that the steady-state gas flow equations may lead to pressure violations in numerical experiments. The gas flow dynamics are formulated by partial differential equations (PDEs) \cite{Dorin_PDEs}. Linear discrete approximation is proposed in \cite{CorreaPosada2015}, in which OGPF is cast as an MILP. \cite{Chaudry_DCOPF} and \cite{Cong_PDE} use coarse approximations to make the PDEs trackable. \cite{Anatoly_PDE} directly use the PDEs to depict gas flow dynamics and solve OGPF using a nonlinear solver.

Recent studies try to replace the Weymouth equation with relaxed but simpler constraints, and leverage the computation superiority of convex optimization. For example, linear relaxation is discussed in \cite{Hantao_LP} and second order cone relaxation is investigated in \cite{Conrado_MISOCP}. Nevertheless, the solution obtained from the relaxed model is not feasible when the convex relaxation is not exact. Provable exactness guarantee is non-trivial.

In the aforementioned studies \cite{Chaudry_DCOPF,Cong_Dispatch,Cong_PDE,Anatoly_PDE,Chengcheng_MILP,Alberto_Newton,Sheng_Probabilistic,Carlos_MILP}, one common assumption is that the power system and gas system are operated in a centralized manner: a central agency that has full control authority over both systems, which may not be consistent with current engineering practice. Furthermore, the power and gas system may be unwilling to exchange their real-time operating data due to privacy and security consideration.  In this regard, a distributed OGPF method with limited data exchange is desired. Pioneer work on distributed OPF problem can be found in \cite{DOPF1,DOPF2,DOPF3}. In \cite{He_ADMM}, a distributed robust optimization framework is designed to solve the day-ahead co-scheduling problem of power and gas system in transmission level, which is related to the topic of this paper, nevertheless, the main focus of this paper is to enhance the computation efficiency through convex optimization, while neglecting uncertainties.

In this paper, a multi-period distribution-level OGPF model is proposed. The power distribution network (PDN) and gas distribution network (GDN) are radial in topology and coupled by gas-fired distributed generators (DGs). The optimal power flow (OPF) of PDN is formulated based on the the branch flow model \cite{Wu_Branch,Wu_Branch2}, which directly incorporates bus voltage magnitude and line power flow. The second order cone (SOC) relaxation \cite{Low_PDN_SOCP} is performed on the non-convex line apparent power equations, and the exactness has been proved in \cite{Na_SOCP_Exact}. In this way, the OPF problem is cast as an second order cone program (SOCP). The optimal gas flow (OGF) of GDN incorporates multi-period and linepack. The gas flow directions are fixed in the dispatch time scale \cite{Nico_SingleDirection} due to the tree topology of GDN. As a result, the sign function in the Weymouth equation can be removed which yields a quadratic equality in the form of a difference of two convex functions. Moreover, each non-convex equality is replaced by a second-order cone inequality and a convex-concave inequality. Following the method in \cite{Boyd_CCP}, OGF can be solved by a sequential SOCP algorithm (SSA) via updating the linear approximation of the concave part. Compared with the state-of-the-art approaches for the Weymouth equation, SSA have two advantages:
\begin{enumerate}
\item Fast and reliable computation. SSA for OGF consists of solving a sequence of SOCPs, which are easier to solve than nonlinear programs or MILPs ;
\item Feasibility guarantee, the feasibility of the solution can be guaranteed by adjusting the convergence criterion;
\end{enumerate}

The contributions of the paper are twofold:

1) SSA for OGF. Convergence and optimality are guaranteed by the theory in \cite{Boyd_CCP}. A dedicated initial point selection scheme is designed, which helps accelerate convergence and improve the solution quality.

2) An alternating direction multiplier method (ADMM) based distributed algorithm for the multi-period distribution-level OGPF problem considering gas dynamics. The OPF and OGF subproblems can be solved independently by the operators of PDN and GDN, while preserving private information of each network.

The rest parts of the paper are structured as follows: Section II introduces the mathematical formulation of OGPF. Section III develops the solution methods. Section IV presents numerical results on two testing systems which validate the proposed algorithm. Finally, Section V concludes the paper.

\section{Mathematical Formulation}
\subsection{Assumptions \& Simplifications}

Before proceeding to mathematical formulation, some prerequisite assumptions and simplifications to facilitate model formulation are stated as follows:

For the modeling of PDN, we assume
\begin{enumerate}
\item PDN is radial. The power flow direction is fixed; reverse power flow is prohibited.
\item The gas flow in the fuel pipeline from GDN to gas-fired DGs of PDN is fully controllable, i.e., the gas demand of gas-fired DGs only depends on its active power output;
\item The electricity demand of PDN is supplied by gas-fired and non-gas DGs.
\end{enumerate}

For the modeling of GDN, we assume
\begin{enumerate}
\item GDN is radial. The gas flow direction is fixed \cite{Nico_SingleDirection};
\item The gas flow dynamics can be approximated by algebraic equations. Details can be found in \cite{CorreaPosada2015};
\item The simplified compressor model in \cite{Wolf_BNG} is used. For detailed compressor model, please refer to \cite{Cong_Dispatch};
\item The cost function of the compressor is linear. According to the discussions in \cite{Wu_3-5} and \cite{Conrado_3-5}, the compressor typically consumes about 3-5\% of the transported gas;
\item The gas demand of GDN is supplied by gas retailers.
\end{enumerate}

\subsection{Centralized OGPF Problem}

The objective of OGPF is to minimize the total production cost of the coupled system, including the fuel cost of non-gas units and gas purchase cost of the GDN. The detailed expression is given in (\ref{central_obj}), where the first component represents the fuel cost of non-gas units, which is usually a convex quadratic function of $p_{gt}$; the second component represents the gas purchase cost. We introduce a constant $\alpha_{c}\in$[0.03,0.05] to quantify the energy consumed by compressor. If it is fueled by gas, then the in/out gas flow of compressor would be different; if it is powered by electricity, then there would be additional variable demand in PDN. Please refer to (\ref{active_balance}), (\ref{gas_balance}), (\ref{com_flow_out_1}) for detailed description.

\begin{equation}\label{central_obj}
    Obj=\min(\sum_{t}(\sum_{g}Q_g(p_{gt})+\sum_{w}Q_wy_{wt}))
\end{equation}
\begin{equation}\label{active_power_capacity}
  P_{min}^{\{\cdot\}}\le p_{\{\cdot\}t} \le P_{max}^{\{\cdot\}},~~\{\cdot\}=\{g,n\},~\forall t
\end{equation}
\begin{equation}\label{reactive_capacity}
  Q_{min}^{\{\cdot\}}\le q_{\{\cdot\}t} \le Q_{max}^{\{\cdot\}},~~\{\cdot\}=\{g,n\},~\forall t
\end{equation}
\begin{equation}\label{voltage}
  (v^{min}_{i_p})^2 \le \nu_{i_pt} \le (v^{max}_{i_p})^2
\end{equation}
\begin{equation}\label{current}
  0 \le I_{l_pt} \le (I_{l_p}^{max})^2
\end{equation}
\begin{equation}
\begin{split}
\label{active_balance}
& \sum_{\{\cdot\}\in \Psi_{\{\cdot\}}(i_p)}p_{\{\cdot\}t}+\sum_{l\in \Psi_{O_2}(i_p)}(pf_{l_pt}-r_{l_p}I_{l_pt})-g_{i_p}\nu_{i_p}\\
&-\sum_{l\in\Psi_{O_1}(i_p)}pf_{l_pt}-\sum_{d\in \Psi_{d_p}(i_p)}p_{d_pt}-\sum_{c\in\Psi_c(i_p)}\chi\alpha_cy_{ct}^{in}\\
&=0,~~\{\cdot\}=\{g,n\},~\forall i_p,t
\end{split}
\end{equation}
\begin{equation}
\begin{split}
\label{reactive_balance}
&\sum_{\{\cdot\}\in \Psi_{\{\cdot\}}(i_p)}q_{\{\cdot\}t}+\sum_{l\in \Psi_{O_2}(i_p)}(qf_{l_pt}-x_{l_p}I_{l_pt})-b_{i_p}\nu_{i_p}-\\
&\sum_{l\in\Psi_{O_1}(i_p)}qf_{l_pt}-\sum_{d\in \Psi_{d_q}(i_p)}q_{d_qt}=0,~\{\cdot\}=\{g,n\},~\forall i_p,t
\end{split}
\end{equation}
\begin{equation}\label{voltage_decreament}
\begin{split}
&\nu_{i_{p2}t}=\nu_{i_{p1}t}-2(r_{l_p}pf_{l_p}+x_{l_p}qf_{l_p})+(r_{l_p}^2+x_{l_p}^2)I_{l_pt},\\
&~\forall (i_{p1},i_{p2})\in l_p, t
\end{split}
\end{equation}
\begin{equation}\label{Apparent_Power}
  I_{l_pt} \ge \frac{pf_{l_pt}^2+qf_{l_pt}^2}{\nu_{i_{p1}t}},~\forall (i_{p1},i_{p2})\in l_p, t
\end{equation}
\begin{equation}\label{gas_output}
  y_w^l \le y_{wt} \le y_w^u~~\forall w,t
\end{equation}
\begin{equation}\label{node_pressure}
\tau_{i_g}^l\le u_{i_gt} \le \tau_{i_g}^u~~\forall i_g,t
\end{equation}
\begin{equation}
\begin{split}
\label{gas_balance}
&\sum_{w\in_{\Theta_w(i_g)}}y_{wt}-\sum_{d_g\in_{\Theta_{d_g}(i_g)}}y_{d_gt}-\sum_{\{\cdot\}\in\Theta_{\{\cdot\}_{1}}(i_g)}y^{out}_{\{\cdot\}t}=\\
&\sum_{n\in_{\Theta_n(i_g)}}p_{nt}/\beta_n-\sum_{\{\cdot\}\in\Theta_{\{\cdot\}_{2}}(i_g)}y^{in}_{\{\cdot\}t}~~\{\cdot\}=\{c,l_g\},~\forall i_g,t
\end{split}
\end{equation}
\begin{equation}\label{weymouth}
\frac{(y_{l_gt}^{in}+y_{l_gt}^{out})|y_{l_gt}^{in}+y_{l_gt}^{out}|}{4}=\phi_{l_g}((u_{l_g^1t})^2-(u_{l_g^2t})^2)~\forall l_g,t
\end{equation}
\begin{equation}\label{wey_co}
\phi_{l_g}=\frac{\pi^2\lambda^2R_{l_g}^5}{16X_{l_g}F_{l_g}\mu T_kZ_{l_g}\rho_0^2}
\end{equation}
\begin{equation}\label{gas_mass}
m_{l_gt}=\frac{\pi}{4}\frac{X_{lg}R_{lg}^2}{\mu T_kZ_{l_g}\rho_0}\frac{u_{l_g^1t}+u_{l_g^2t}}{2}~~\forall l_g,t
\end{equation}
\begin{equation}\label{gas_mass_storage}
m_{l_gt}=m_{l_g(t-1)}+y^{in}_{l_gt}-y^{out}_{l_gt}~~\forall l_g,t
\end{equation}
\begin{equation}\label{com_pressure}
  u_{i_{gt,c_2}}\le \gamma_cu_{i_{gt,c_1}}~~\forall c,t
\end{equation}
\begin{equation}\label{com_flow_in}
  0\le y_{ct}^{in} \le y_c^{max}~~\forall c,t
\end{equation}
\begin{equation}\label{com_flow_out_1}
  y_{ct}^{in}-(1-\alpha_c)y_{ct}^{out}=0~~\forall c,t
\end{equation}

\vspace{6pt}
\noindent Constraint (\ref{active_power_capacity})-(\ref{com_flow_out_1}) represent the operating constraints of the interconnected system. For the PDN, (\ref{active_power_capacity}) and (\ref{reactive_capacity}) enforce the capacity of active and reactive generation. (\ref{voltage}) and (\ref{current}) restrict the boundary of nodal voltage square and line current square. (\ref{active_balance}) and (\ref{reactive_balance}) are nodal power balancing conditions. (\ref{voltage_decreament}) is the voltage drop equation. (\ref{Apparent_Power}) expresses the relaxed line apparent power. It is proved that (\ref{Apparent_Power}) will be active at the optimal solution  \cite{Low_PDN_SOCP}, thus the relaxation is exact. For the GDN, (\ref{gas_output}) and (\ref{node_pressure}) restrict the gas supply range of the retailer and gas pressure at each node. (\ref{gas_balance}) is the nodal gas balancing equation, where $\Theta_w(i_g)$ ($\Theta_{d_g}(i_g)$, $\Theta_{n}(i_g)$) represent the set of gas retailer (gas load, gas-fired DGs) connecting to node $i_g$; $\Theta_{c_{1}}(i_g)$, $\Theta_{c_{2}}(i_g)$, $\Theta_{l_{g_{1}}}(i_g)$, $\Theta_{l_{g_{2}}}(i_g)$ represent the set of active/passive pipelines whose head/tail node is $i_g$. Particularly, pipelines with and without compressors are referred to as active and passive pipelines, respectively. (\ref{weymouth}) is the Weymouth equation which characterizes the relationship between gas flow in a passive pipeline and node pressure, where $u_{l_g^1t}$, $u_{l_g^2t}$ are the pressure of initial and terminal node of $l_g$, respectively. (\ref{wey_co}) defines the coefficient $\phi_{l_g}$ in the Weymouth equation. (\ref{gas_mass}) gives the relationship between line pack and average pressure of the line. (\ref{gas_mass_storage}) depicts the time-dependent relationship between line pack and gas flow. The mathematical procedure to derive (\ref{gas_mass}) and (\ref{gas_mass_storage}) can be found in the Appendix of \cite{CorreaPosada2015}, in which a constant compressibility factor $Z_{l_g}$ is adopted to preserve linearity. For each active pipeline, (\ref{com_pressure}) and (\ref{com_flow_in}) limits the maximum compression ratio and gas flow. (\ref{com_flow_out_1}) indicates the relationship between in/out gas flow of an active pipeline. In fact, (\ref{com_pressure})-(\ref{com_flow_out_1}) are simplified model for compressors. The exact model of compressor in \cite{Cong_Dispatch} is highly nonlinear and non-convex.

The sign function in (\ref{weymouth}) greatly challenges the solution. According to the gas system operation practice, the gas flow direction in GDN does not change intra-day \cite{Nico_SingleDirection}. In this regard, (\ref{weymouth}) can be reduced as
\begin{equation}\label{weymouth_easy}
(y_{l_gt}^{in}+y_{l_gt}^{out})^2/4 = \phi_{l_g}((u_{l_g^1t})^2-(u_{l_g^2t})^2),~~\forall l_g,t
\end{equation}
\begin{equation}\label{weymouth_easy_positive}
u_{l_g^1t} \ge u_{l_g^2t} \ge 0, ~~\forall l_g,~t
\end{equation}

\noindent where we assume the notation of initial and terminal node of $l_g$ is consistent with the positive direction of gas flow. (\ref{weymouth_easy}) and (\ref{weymouth_easy_positive}) hold naturally for tree-like GDNs. The non-convexity in OGF originates from equation (\ref{weymouth_easy}).

\subsection{Distributed Counterpart of OGPF Problem}
In practice, it may not be desired to solve the OGPF in a centralized manner for several reasons. On the one hand, the GDN and PDN are operated by different entities. There is no central agency which has full control authority of both infrastructures. On the other hand, the operators of GDN and PDN may not wish to share their private information for security consideration. In the proposed framework, OGPF is decomposed into an OPF subproblem and an OGF subproblem, which will be solved independently.

\subsubsection{OPF subproblem}
The operator of PDN endeavours to minimize the fuel cost of non-gas DGs subject to PDN constraints with fixed demands from compressors, rendering
\begin{equation}\label{Obj_OPF}
Obj_{power}=\min \sum_{t}\sum_{g}Q_g(p_{gt})
\end{equation}
\begin{equation*}
~~~s.t.~(\ref{active_power_capacity})-(\ref{current}), (\ref{reactive_balance})-(\ref{Apparent_Power})
\end{equation*}
\begin{equation}
\begin{split}
\label{active_balance_power}
& \sum_{\{\cdot\}\in \Psi_{\{\cdot\}}(i_p)}p_{\{\cdot\}t}+\sum_{l\in \Psi_{O_2}(i_p)}(pf_{l_pt}-r_{l_p}I_{l_pt})-g_{i_p}\nu_{i_p}-\\
&\sum_{l\in\Psi_{O_1}(i_p)}pf_{l_pt}-\sum_{d\in \Psi_{d_p}(i_p)}p_{d_pt}=0,~\{\cdot\}=\{g,n\},~\forall i_p,t
\end{split}
\end{equation}
(\ref{active_balance_power}) is the active power balance equation of PDN nodes, which are not connected with compressors of GDN. Compared with (\ref{active_balance}), electricity consumed by compressors are not involved.

\subsubsection{OGF subproblem}
The operator of GDN seeks to minimize the gas purchase cost subject to GDN constraints with fixed demands from gas-fired units, rendering
\begin{equation}\label{Obj_OGF}
Obj_{gas}=\min\sum_{t}\sum_{w}Q_wq_{wt}
\end{equation}
\begin{equation*}
~~~s.t.~(\ref{gas_output})-(\ref{node_pressure}), (\ref{wey_co})-(\ref{weymouth_easy_positive})
\end{equation*}
\begin{equation}
\label{gas_balance_gas}
\begin{split}
&\sum_{w\in_{\Theta_w(i_g)}}y_{wt}-\sum_{d_g\in_{\Theta_{d_g}(i_g)}}y_{d_gt}=\sum_{\{\cdot\}\in\Theta_{\{\cdot\}_{1}}(i_g)}y^{out}_{\{\cdot\}t}\\
&-\sum_{\{\cdot\}\in\Theta_{\{\cdot\}_{2}}(i_g)}y^{in}_{\{\cdot\}t}~~\{\cdot\}=\{c,l_g\},~\forall i_g,t
\end{split}
\end{equation}
(\ref{gas_balance_gas}) is the gas balance equation of GDN nodes, which are not connected with gas-fired DGs of PDN. Compared with (\ref{gas_balance}), gas consumed by gas-fired DGs in PDN are not involved.

\subsubsection{Coupling Constraints}
As the PDNs and GDNs are interconnected by gas-fired units and electrical compressors, the following coupling constraints should be satisfied.
\begin{equation}
\label{active_balance_couple}
\begin{split}
&\sum_{\{\cdot\}\in \Psi_{\{\cdot\}}(i_p)}p_{\{\cdot\}t}+\sum_{l\in \Psi_{O_2}(i_p)}(pf_{l_pt}-r_{l_p}I_{l_pt})-g_{i_p}\nu_{i_p}\\
&-\sum_{l\in\Psi_{O_1}(i_p)}pf_{l_pt}-\sum_{d\in \Psi_{d_p}(i_p)}p_{d_pt}-\sum_{c\in\Psi_c(i_p)}\chi_c\alpha_cy_{ct}^{in}\\
&=0,~~\{\cdot\}=\{g,n\},~\Psi_c(i_p)\neq\emptyset,~\forall t
\end{split}
\end{equation}
\begin{equation}
\label{gas_balance_couple}
\begin{split}
&\sum_{w\in_{\Theta_w(i_g)}}y_{wt}-\sum_{d_g\in_{\Theta_{d_g}(i_g)}}y_{d_gt}-\sum_{\{\cdot\}\in\Theta_{\{\cdot\}_{1}}(i_g)}y^{out}_{\{\cdot\}t}=\\
&\sum_{n\in_{\Theta_n(i_g)}}p_{nt}/\beta_n-\sum_{\{\cdot\}\in\Theta_{\{\cdot\}_{2}}(i_g)}y^{in}_{\{\cdot\}t},~~\{\cdot\}=\{c,l_g\},\\
&~~\Theta_n(i_g)\neq\emptyset,~\forall t
\end{split}
\end{equation}
(\ref{active_balance_couple}) and (\ref{gas_balance_couple}) are active power and gas balance equations of PDN and GDN nodes, respectively, which are connected with compressors and gas-fired DGs, respectively.

\section{Solution Methodology}
In this section, a distributed algorithm is proposed to solve OGPF iteratively. In particular, a sequential SOCP based method is devised to solve the non-convex OGF subproblem. For notation simplicity, the compact form of the centralized OGPF problem is written as
\begin{equation}\label{compact_form}
\begin{split}
&\min F_p(\vec{x})+F_g(\vec{z})\\
s.t.~~&\vec{x}\in \Omega_{x},~\vec{z}\in\Omega_{z}\\
& \vec{Ax}+\vec{Bz}=\vec{c}
\end{split}
\end{equation}
In (\ref{compact_form}), $\vec{x}$ and $\vec{z}$ are the decision variables of the OPF subproblem and the OGF subproblem, respectively. $F_p(\cdot)$ and $F_g(\cdot)$ are the corresponding objective functions. $\Omega_x$ and $\Omega_z$ are the feasible region composed by the constraints of OPF and OGF part, respectively. The last equality constraint represents the coupling constraints in (\ref{active_balance_couple}) and (\ref{gas_balance_couple}), where matrices $\vec{A,B,c}$ can be extracted from the respective coefficients of variables. It can be observed that the only nonconvexity in (\ref{compact_form}) stems from the nonlinear Weymouth equation (\ref{weymouth_easy}) in $\Omega_z$. Next we will first present the distributed algorithm for OGPF, then derive a method to solve the non-convex OGF subproblem.

\subsection{Distributed Algorithm for OGPF}
The objective function of (\ref{compact_form}) can be decomposed into convex functions  $F_p$ and $F_g$. Feasible regions defined by $\Omega_x$ and coupling constraints are convex. Most constraints in $\Omega_z$ are also convex. In view of these observations, the ADMM \cite{Boyd_ADMM} is adopted to solve OGPF in a distributed manner. The augmented Lagrangian function of (\ref{compact_form}) is expressed as
\begin{equation}\label{Augmented_Lagrangian}
\begin{split}
  L(\vec{x},\vec{z},\vec{\xi})=&F_p(\vec{x})+F_g(\vec{z})+\vec{\xi}^T(\vec{Ax}+\vec{Bz}-\vec{c})\\
  &+\frac{d}{2}||\vec{Ax}+\vec{Bz}-\vec{c}||^2
\end{split}
\end{equation}
where $\vec{\xi}$ is the dual variable of the coupling constraints in (\ref{compact_form}) and $d$ is a constant parameter. Details are given in Algorithm~1.

\begin{algorithm}[!t]
\caption{ADMM for OGPF}
\label{alg:ADMM}
\begin{algorithmic}[1]
\STATE Given initial values of $\vec{x}_0$ and $\vec{\xi}_0$. Set iteration index $k=0$, select value for $k_{max}, $constant $d$ and $\varsigma$.
\STATE Solve problems (\ref{update_x}) and (\ref{update_z})
\begin{equation}\label{update_x}
  \vec{x}_{k+1}=\arg{\min_{\vec{x}\in\Omega_x}(L(\vec{x},\vec{z}_k,\vec{\xi}_k)}-F_g(\vec{z}_k))
\end{equation}
\begin{equation}\label{update_z}
  \vec{z}_{k+1}=\arg{\min_{\vec{z}\in\Omega_z}(L(\vec{x}_{k+1},\vec{z},\vec{\xi}_k)}-F_p(\vec{x}_{k+1}))
\end{equation}
update $\vec{\xi}$ according to (\ref{update_xi})
\begin{equation}\label{update_xi}
  \vec{\xi}_{k+1}=\vec{\xi}_k+d(\vec{Ax}_{k+1}+\vec{Bz}_{k+1}-\vec{c})
\end{equation}

\STATE Set $k=k+1$. If $\max|\vec{Ax}_{k+1}+\vec{Bz}_{k+1}-\vec{c}|\le\varsigma$ is satisfied, then quit; else if $k=k_{max}$, then quit; else, go to Step 2.
\end{algorithmic}
\end{algorithm}

The convergence guarantee of ADMM relies on the model convexity \cite{Boyd_ADMM}. However, the non-convex Weymouth equation challenges the convergence of Algorithm \ref{alg:ADMM}. In view of its quadratic convex-concave structure, a sequential SOCP method is devised to solve the OGF subproblem in the next subsection, such that and the global optimal solution can be efficiently computed. The local convexity of OGF also contributes to the convergence of Algorithm 1.

\subsection{Convexification of OGF Part}

The non-convex Weymouth equation is a set of quadratic equality in the form of the difference of convex functions, which can be cast as opposite inequalities as
\begin{equation}\label{weymouth_easy_re1}
(y_{l_gt}^{in}+y_{l_gt}^{out})^2/4 + \phi_{l_g}(u_{l_g^2t})^2 - \phi_{l_g}(u_{l_g^1t})^2 \le 0~~\forall l_g,t
\end{equation}
\begin{equation}\label{weymouth_easy_re2}
  \underbrace{\phi_{l_g}(u_{l_g^1t})^2}_{f_{l_gt}}-\underbrace{((y_{l_gt}^{in}+y_{l_gt}^{out})^2/4+\phi_{l_g}(u_{l_g^2t})^2}_{g_{l_gt}})\le 0~~\forall l_g,t
\end{equation}
where (\ref{weymouth_easy_re1}) is in fact a set of SOC inequalities whose canonical form is given by
\begin{equation}
\label{weymouth_socp}
\left\|  \begin{array}{*{20}c}
   (y_{l_gt}^{in}+y_{l_gt}^{out})/2  \\
   \sqrt{\phi_{l_g}} u_{l_g^2t}  \\
\end{array}   \right\|  \le  \sqrt{\phi_{l_g}} u_{l_g^1t}
\end{equation}

Given a vector $[(y_{l_gt,k}^{in})^T~(y_{l_gt,k}^{out})^T~u_{i_gt,k}^T]^T$, the linear approximation of $g_{l_gt}$ in (\ref{weymouth_easy_re2}) is
\begin{equation}
\begin{split}
\label{wey_easy_re2g}
  \hat{g}_{l_gt}&=\frac{(y_{l_gt,k}^{in}+y_{l_gt,k}^{out})(y_{l_gt}^{in}+y_{l_gt}^{out})}{2}-\frac{(y_{l_gt,k}^{in}+y_{l_gt,k}^{out})^2}{4}\\
  &-\phi_{l_g}(u_{l_g^2t,k})^2+2\phi_{l_g}u_{l_g^2t,k}u_{l_g^2t}
\end{split}
\end{equation}

SSA for OGF subproblem is presented in Algorithm 2.
\begin{algorithm}[ht]
\caption{SSA for the OGF Part}
\label{alg:Improved_CCP}
\begin{algorithmic}[1]
\STATE Given $\vec{x}_{k+1}$ and $\vec{\xi}_k$ of $k_{th}$ iteration as well as constant $d$ of Algorithm \ref{alg:ADMM}. Given a initial value of $\vec{z}_{k+1}^{0}$. Set iteration index $j=0$, select value of $\delta, \varrho_0, \varrho_{max}, \epsilon, j^{max}, Obj_{gas}^{0}$ and $\kappa>1$.
\STATE Call (\ref{wey_easy_re2g}). Set the value of $\vec{z}_{k+1}^{j+1},\varsigma_{l_gt}^{+,j+1},\varsigma_{l_gt}^{-,j+1}$ as the solution of the following problem
\begin{spacing}{0.4}
\begin{equation}
\begin{split}\label{Improved_Penalty_DC_basic_con}
&Obj_{gas}^{j+1}=\min_{\vec{z}_{k+1}^{j+1},\varsigma_{l_gt}^{+,j+1},\varsigma_{l_gt}^{-,j+1}}(L(\vec{x}_{k+1},\vec{z}_{k+1}^{j+1},\vec{\xi}_{k})\\
&-F_p(\vec{x}_{k+1})+\varrho^j\sum_{t}\sum_{l_g}\varsigma_{l_gt}^{j+1})\\
s.t.~&f_{l_gt}(\vec{z}_{k+1}^{j+1})-\hat{g}_{l_gt}(\vec{z}_{k+1}^{j+1};\vec{z}_{k+1}^{j})\le \varsigma_{l_gt}^{j+1}~~\forall l_g,t\\
&\varsigma_{l_gt}^{j+1}\ge 0~~\forall l_g,t\\
&(\ref{gas_output})-(\ref{node_pressure}), (\ref{weymouth})-(\ref{com_flow_out_1}), (\ref{gas_balance_gas}),(\ref{weymouth_socp})
\end{split}
\end{equation}
\end{spacing}
\STATE Set $k=k+1$. $\varrho^{j+1}=\min(\kappa\varrho^j,\varrho_{max})$. If both (\ref{Penalty_Stop_DC}) and (\ref{feasibility_criterion}) are satisfied, then quit; else if $j=j^{max}$, then quit; else, go to Step 2.
\begin{spacing}{0.5}
\begin{equation}\label{Penalty_Stop_DC}
  |Obj_{gas}^{j+1}-Obj_{gas}^{j}|\le\delta
\end{equation}

\begin{equation}\label{feasibility_criterion}
\max\{\frac{\varsigma^{j+1}_{l_gt}}{|f_{l_gt}-\hat{g}_{l_gt}|}\}\le\epsilon,~\forall l_g,t
\end{equation}
\end{spacing}
\end{algorithmic}
\end{algorithm}

The convergence proof of Algorithm \ref{alg:Improved_CCP} can be found in \cite{Boyd_CCP}. According to \cite{Boyd_CCP}, (\ref{Penalty_Stop_DC}) can always be satisfied given a large $j^{max}$. However, (\ref{feasibility_criterion}) may not be satisfied even if $j^{max}$ is sufficient large. The reasons are twofold:
\begin{enumerate}
\item Infeasibility of the original problem, which indicates a secure operating condition does not exist. Long-term solution could be upgrading system components, such as investing new pipelines and compressors.

\item Inappropriate initial point selection. Linear approximation performed in (\ref{wey_easy_re2g}) relies on an initial point, which will have a notable impact on the performance of Algorithm 2, such as the number of iterations and solution quality. Zero initiation is the default setting for many heuristic algorithms but may not be a good choice for OGF because the gradient of $g_{-,l_gt}$
is also zero.
\end{enumerate}

\subsection{Initial Point Selection}
An SOCP relaxation based method is provided to obtained a good initial point. Replace (\ref{weymouth_easy}) with the following inequalities

\begin{spacing}{0.5}
\begin{equation}\label{SOCP_wey}
  (y^{in}_{l_gt}+y^{out}_{l_gt})^2/4\le \phi_{l_g}(\zeta_{l_g^1t}-\zeta_{l_g^2t})~~\forall l_g,t
\end{equation}

\begin{equation}\label{SOCP_pressure_1}
  (u_{i_g})^2 \le \zeta_{i_gt}~~\forall i_g,t
\end{equation}

\begin{equation}\label{SOCP_pressure_2}
  (\tau_{i_g}^l)^2\le \zeta_{i_gt} \le (\tau_{i_g}^u)^2~~\forall i_g,t
\end{equation}
\end{spacing}
\vspace{6pt}
\noindent Where $\zeta$ represents the node pressure square. If both (\ref{SOCP_wey}) and (\ref{SOCP_pressure_1}) are equalities, (\ref{SOCP_wey})-(\ref{SOCP_pressure_2}) and (\ref{weymouth_easy}) are equivalent. The initial point of OGF can be obtained by solving (\ref{Augmented_Lagrangian}) with (\ref{gas_output})-(\ref{node_pressure}), (\ref{wey_co})-(\ref{com_flow_out_1}) and (\ref{SOCP_wey})-(\ref{SOCP_pressure_2}), which is a SOCP and can be efficient solved by the off-the-shelf solvers. It should be noted the obtained initial point is not a feasible point, in general, to the OGF problem in Section II.C.

\section{Illustrative Example}
In this section, we present numerical experiments on two test systems to show the effectiveness of the proposed method. The experiments are performed on a laptop with Intel(R) Core(TM) 2 Duo 2.2 GHz CPU and 4 GB memory. The proposed algorithms are coded in MATLAB with YALMIP toolbox. SOCPs are solved by Gurobi 6.5. OGPF in nonlinear form is solved by IPOPT of OPTI Toolbox. The parameters of the solvers are as default without particular mention.

\subsection{13-Bus Power Feeder with 7-Node Gas System}
\begin{figure}
\centering
  \includegraphics[width=0.35\textwidth]{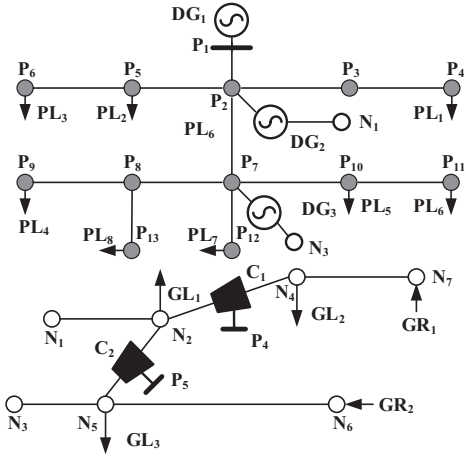}
  \caption{Topology of Power13Gas7 System.}
  \label{fig:Topology_P13G7}
\end{figure}

\begin{figure}
\centering
  \includegraphics[width=0.45\textwidth]{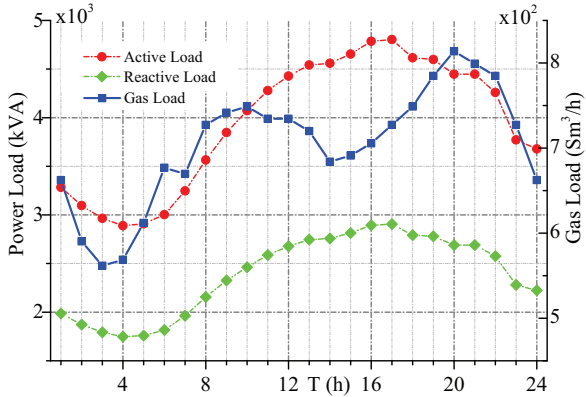}
  \caption{Power and gas demand for Power13Gas7 system.}
  \label{fig:Demand}
\end{figure}

Fig. \ref{fig:Topology_P13G7} depicts the topology of the connected infrastructure. It has 2 gas-fired DGs, 1 non-gas DG, 2 gas retailers, 2 compressors, 4 passive pipelines, 8 power loads and 3 gas loads. In Fig. \ref{fig:Topology_P13G7}, we use $P,PL,DG$ with subscripts to denote the power nodes/power loads/units, respectively, and $N,C,GR,GL$ with subscripts to denote the gas nodes/compressors/gas retailers/gas loads, respectively. Specially, the fuel of $G_2$ and $G_3$ comes from $N_3$ and $N_1$, respectively. Both $C_1$ and $C_2$ are driven by electricity and served by $P_4$ and $P_5$, respectively. The power and gas demand profiles are shown in Fig. \ref{fig:Demand}. The parameters of the system can be found in \cite{Power6Gas7}. This system is referred to as Power13Gas7 system.

\subsection{Simulation Results}
\begin{figure}
\begin{minipage}[t]{0.5\linewidth}
  \centering
  \includegraphics[width=1\textwidth]{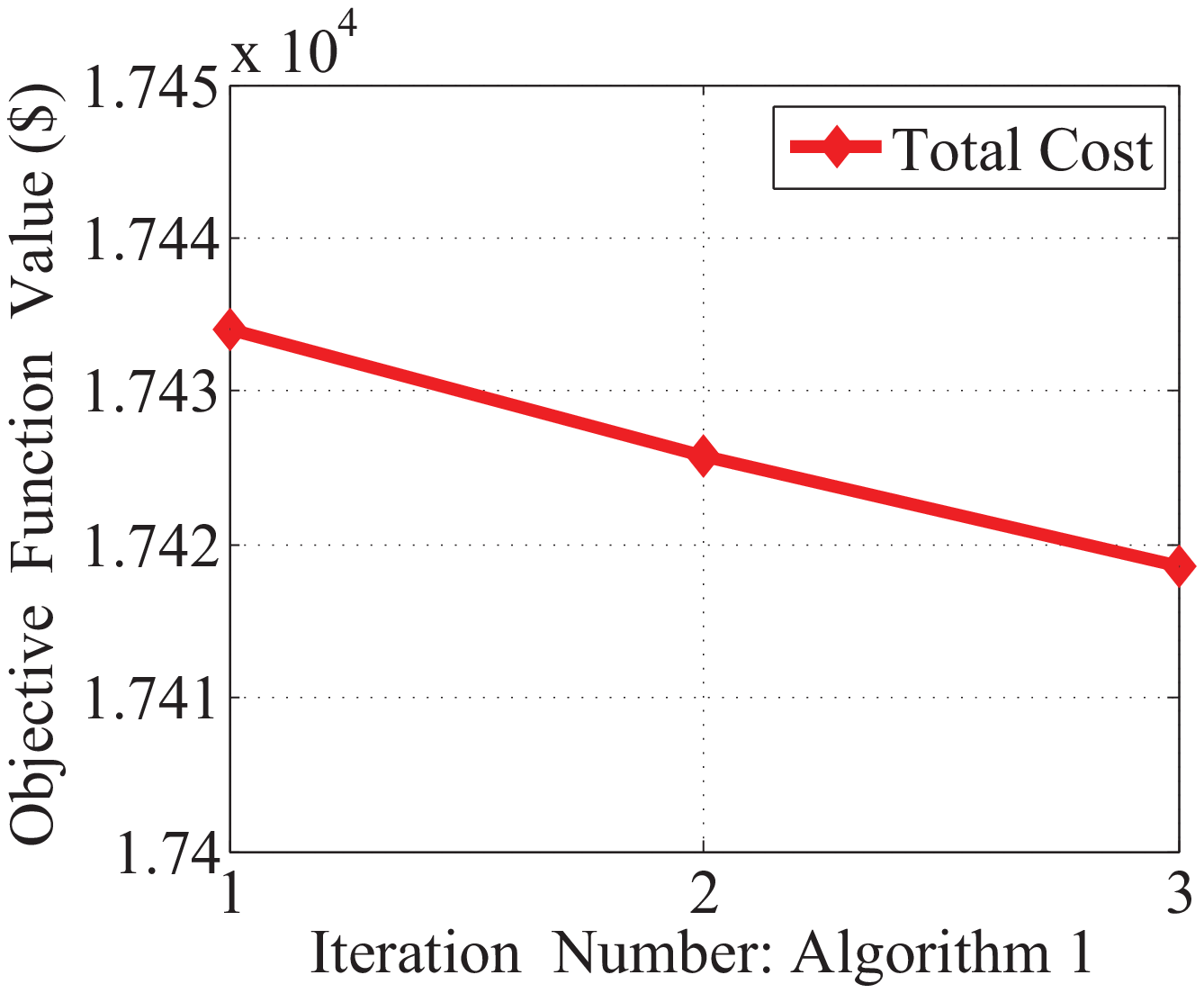}
  \end{minipage}%
  \begin{minipage}[t]{0.5\linewidth}
    \centering
  \includegraphics[width=1\textwidth]{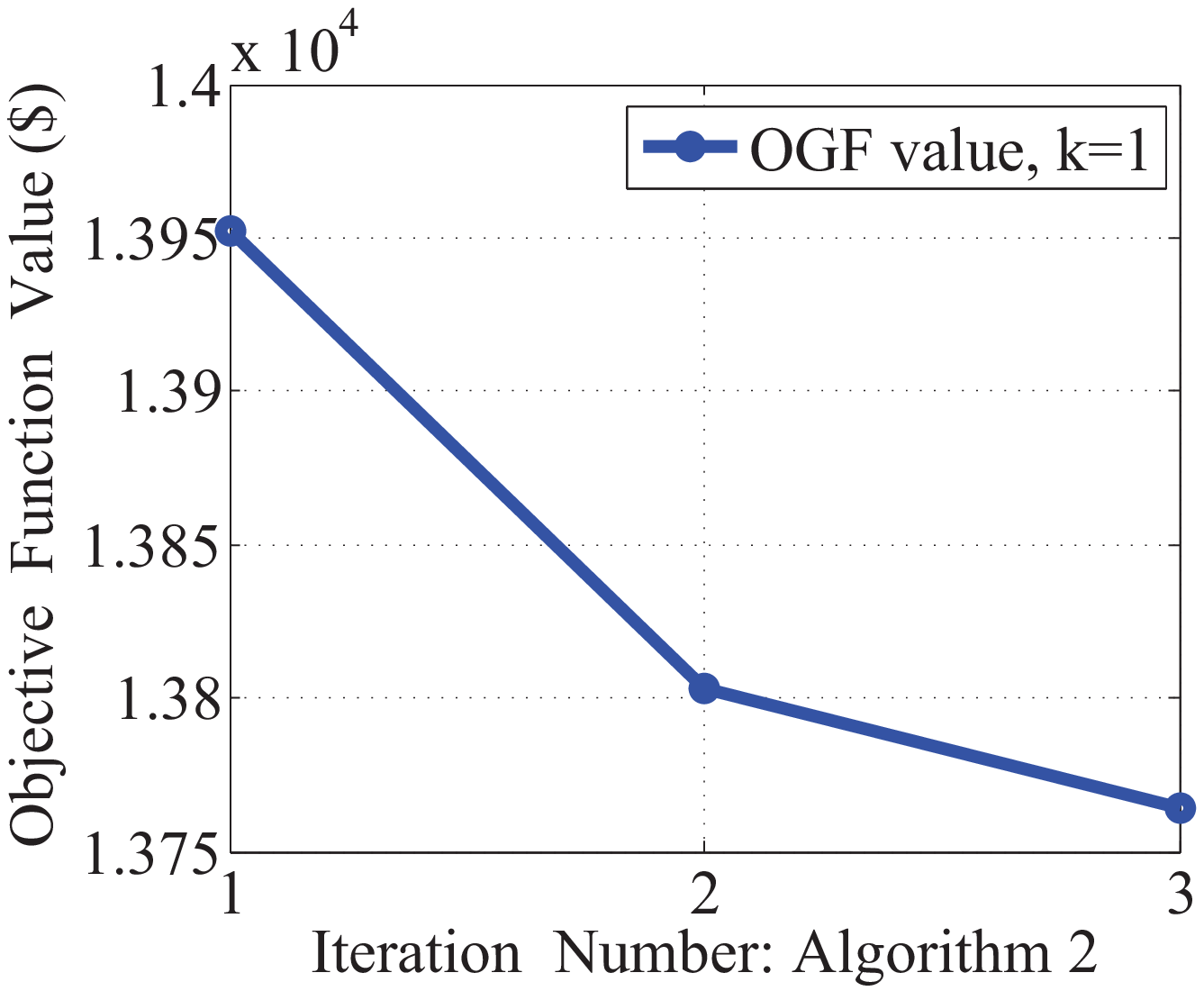}
  \end{minipage}\\
  \begin{minipage}[t]{0.5\linewidth}
    \includegraphics[width=1\textwidth]{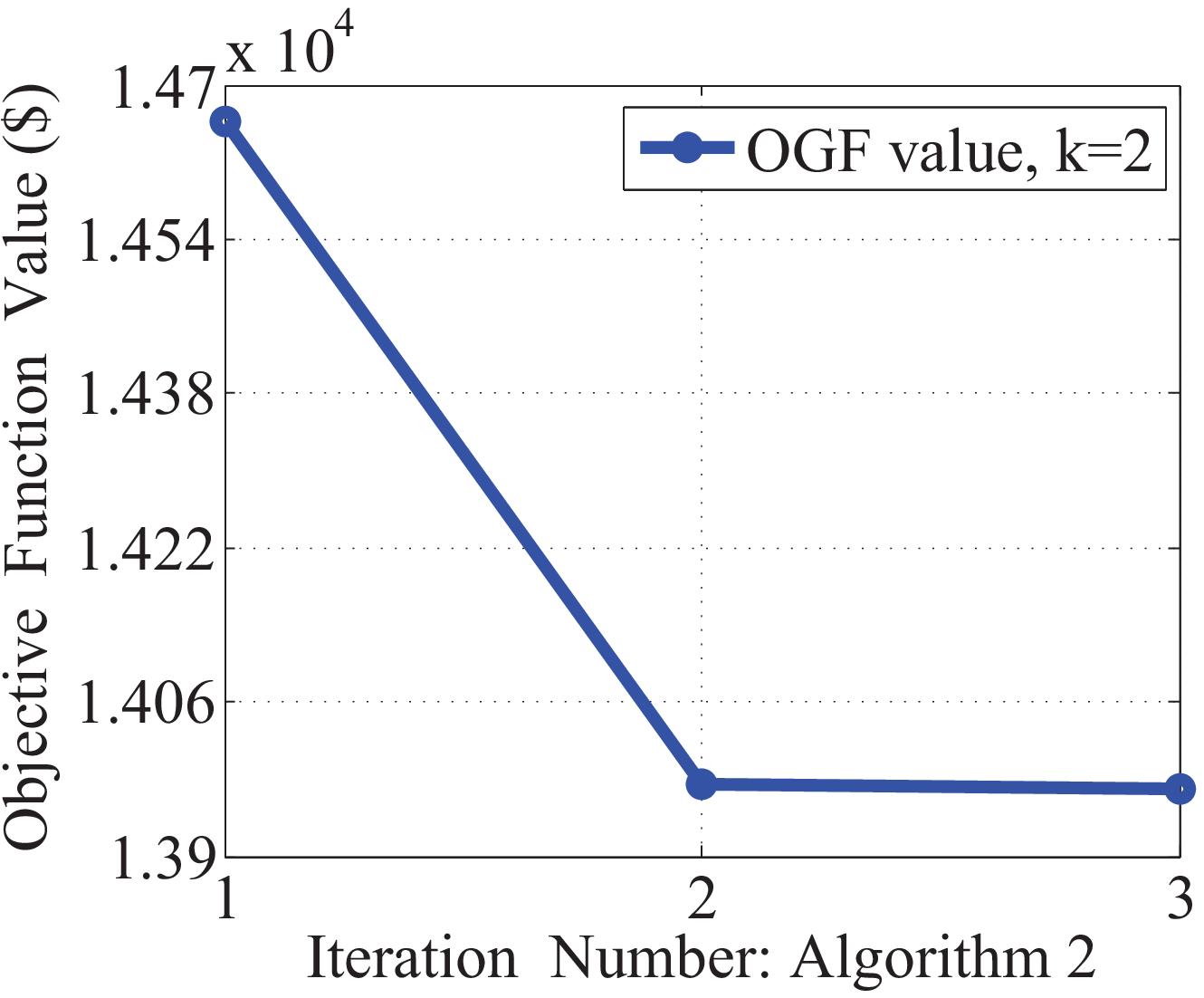}
  \end{minipage}%
  \begin{minipage}[t]{0.5\linewidth}
    \centering
  \includegraphics[width=1\textwidth]{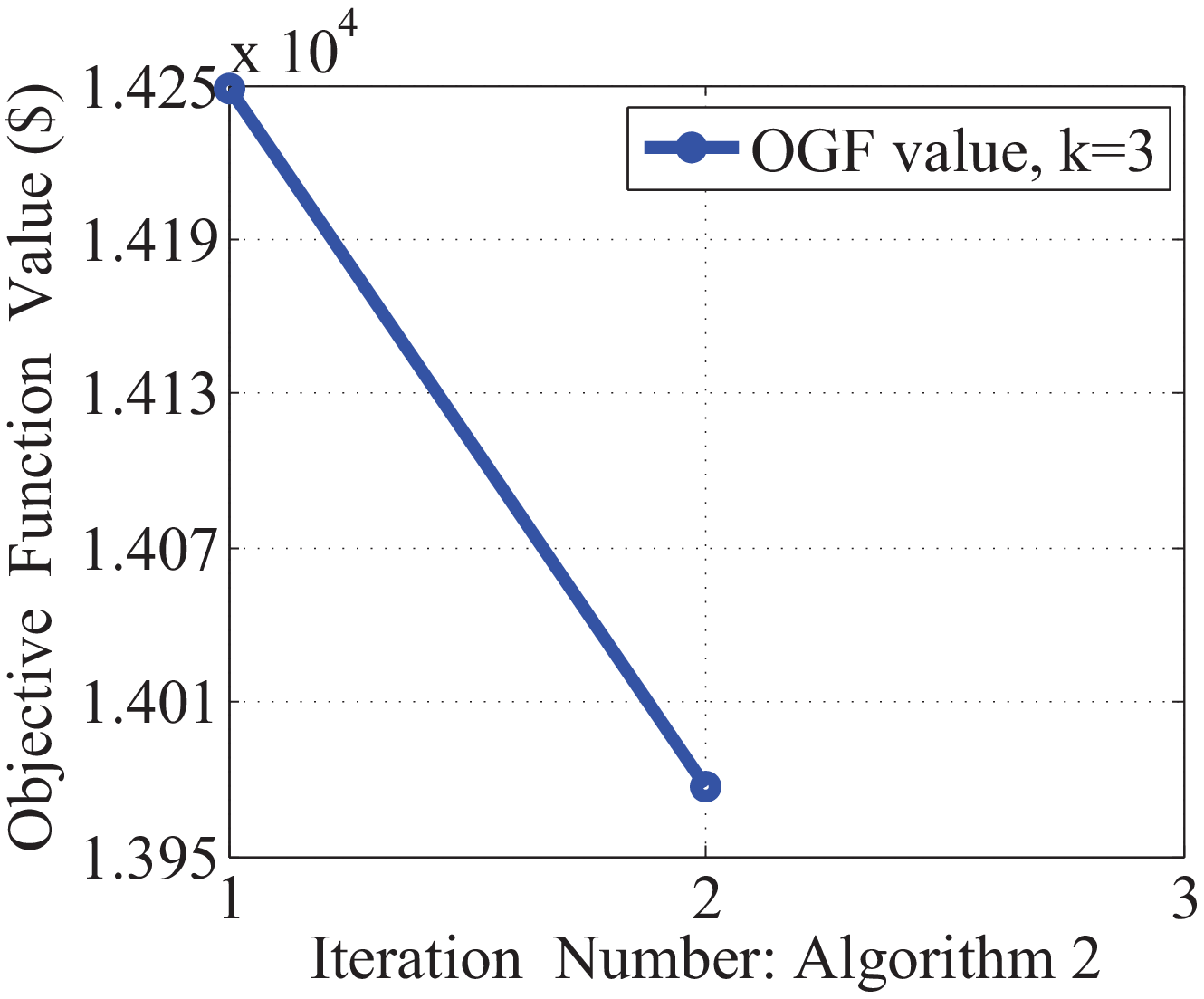}
  \end{minipage}%
  \caption{Optimal value of objective function during iteration}\label{fig:Optimality}
\end{figure}

\begin{figure}
\begin{minipage}[t]{0.5\linewidth}
  \centering
  \includegraphics[width=1\textwidth]{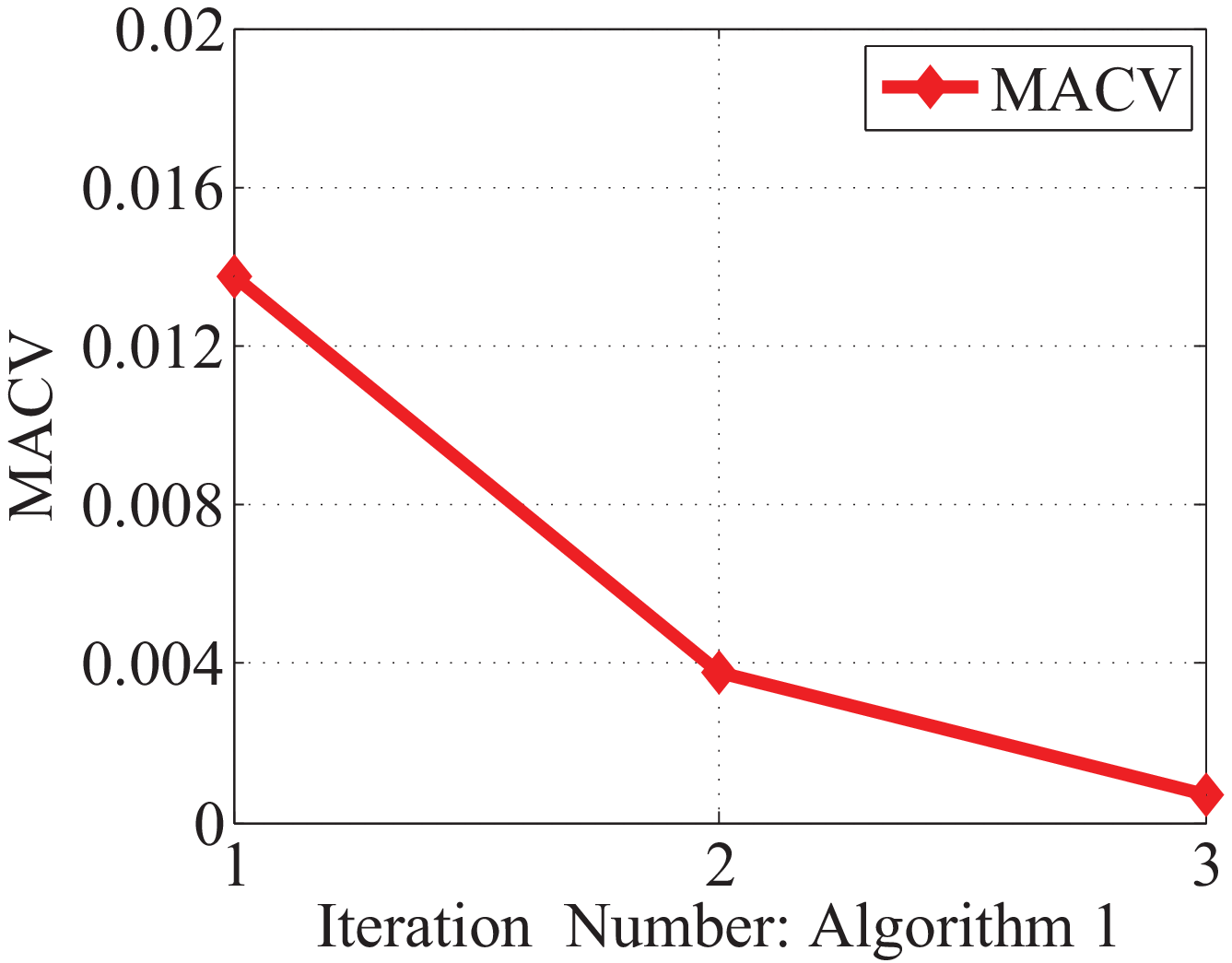}
  \end{minipage}%
  \begin{minipage}[t]{0.5\linewidth}
    \centering
  \includegraphics[width=1\textwidth]{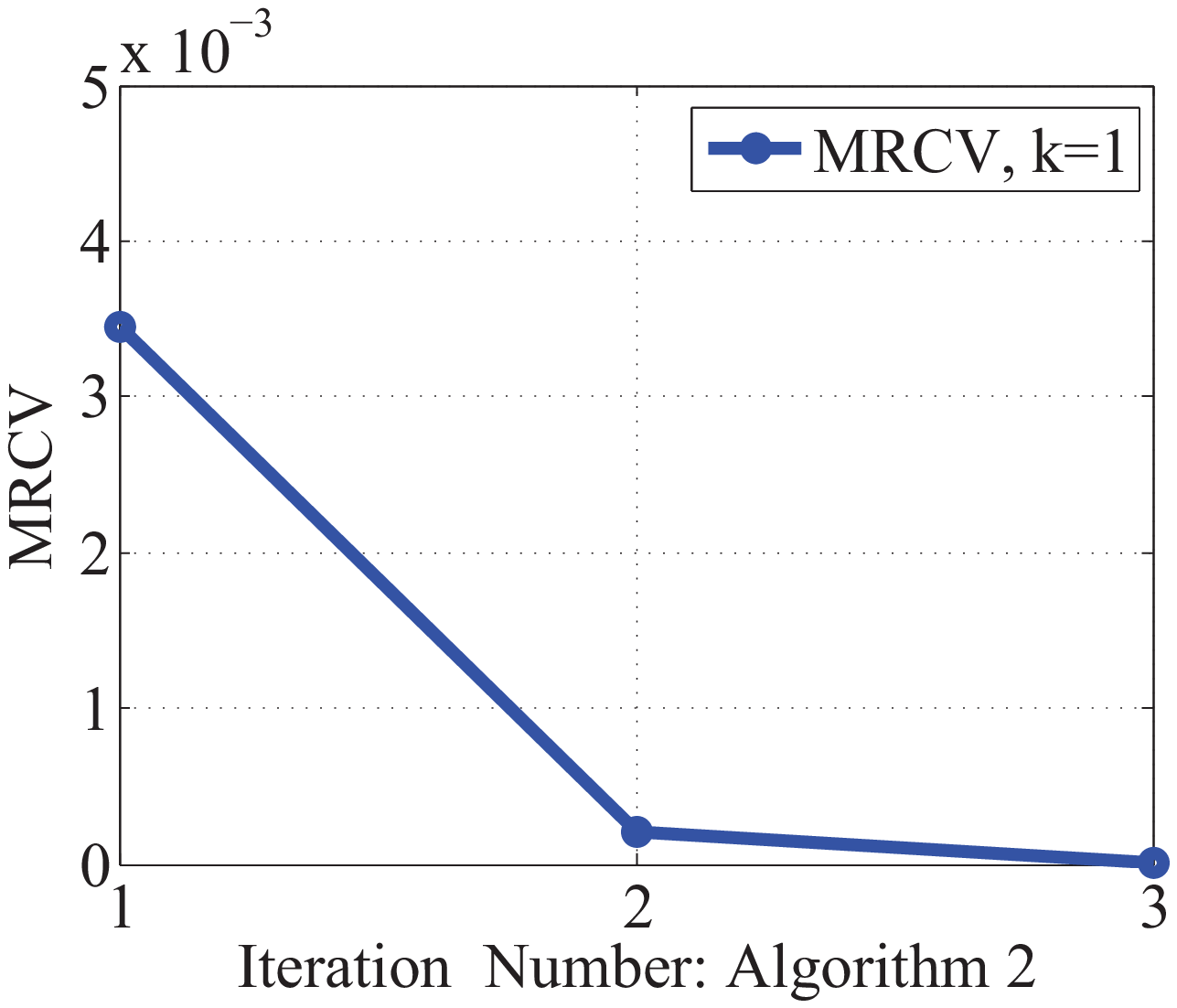}
  \end{minipage}\\
  \begin{minipage}[t]{0.5\linewidth}
    \includegraphics[width=1\textwidth]{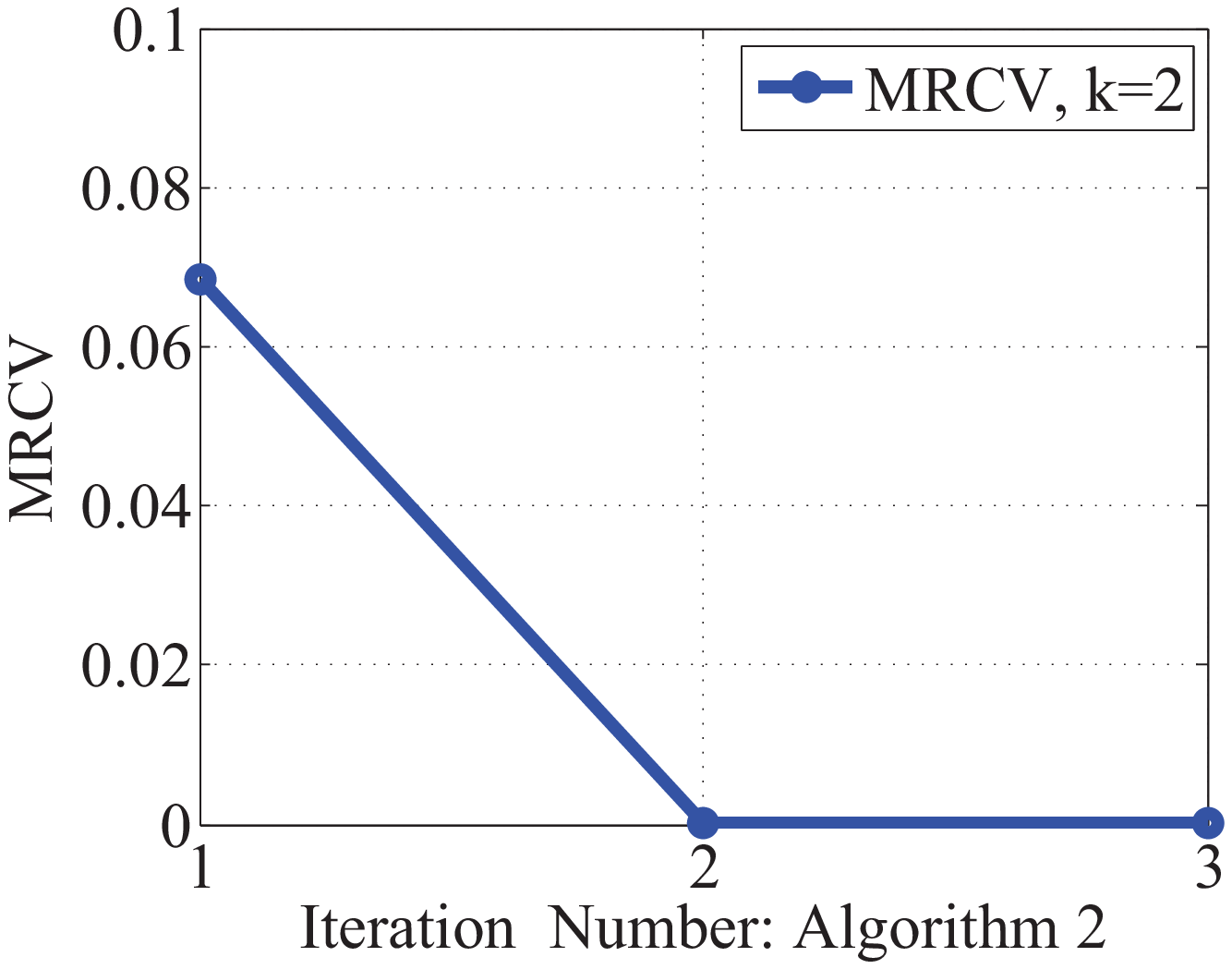}
  \end{minipage}%
  \begin{minipage}[t]{0.5\linewidth}
    \centering
  \includegraphics[width=1\textwidth]{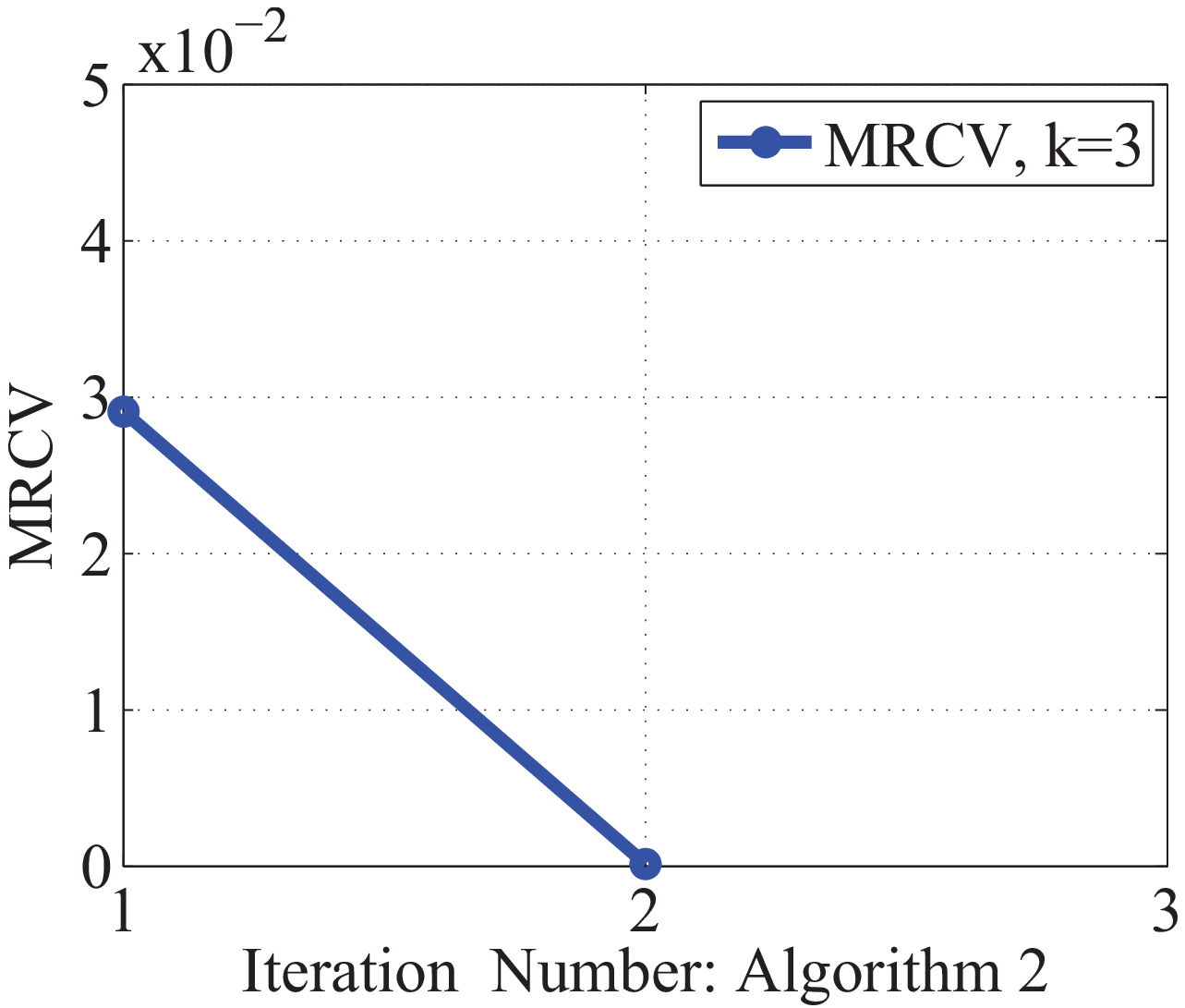}
  \end{minipage}%
  \caption{Feasibility of solution during iteration}\label{fig:Feasibility}
\end{figure}

The start-up parameters of the proposed algorithms are introduced below. The parameters of Algorithm \ref{alg:ADMM} and \ref{alg:Improved_CCP} can be found in Table \ref{Tab:parameter}. The initial value $\vec{z}^{0}_{k+1}$ in Algorithm \ref{alg:Improved_CCP} is obtained by the SOCP-relaxation method proposed in Section III.B. Algorithms \ref{alg:ADMM} and \ref{alg:Improved_CCP} converge in 3 iterations. The objective value sequence generated by  these  algorithms are shown in Fig. \ref{fig:Optimality}. From Fig. \ref{fig:Optimality} we can see that both Algorithm \ref{alg:ADMM} and \ref{alg:Improved_CCP} converge quickly. It should be noted that the total cost of the proposed distributed OGPF is the same with the centralized OGPF model, which is $1.7436\times 10^{4}\$$. From the feasibility perspective, in each iteration of Algorithm \ref{alg:ADMM}, the solution of OGF part strictly satisfies the Weymouth equation, and the final solution offered by Algorithm \ref{alg:ADMM} also satisfies the coupling constraints, as shown in Fig. \ref{fig:Feasibility}. Particularly, MACV and MRCV represent maximum absolute and relative constraint violation, respectively.

\begin{figure}
  \centering
  \includegraphics[width=0.48\textwidth]{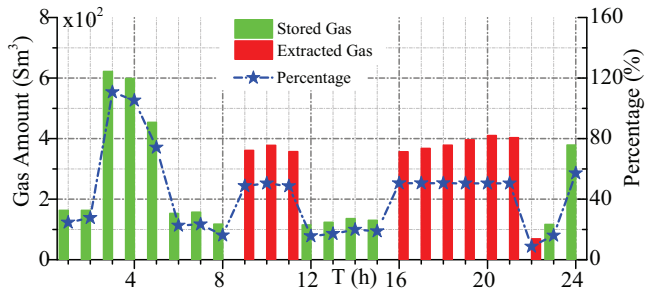}
  \caption{Variation of the sum of line pack during 24 hours}\label{fig:linepack}
\end{figure}

Further, the variation of the sum of line pack in the GDN of Power13Gas7 system are shown in Fig. \ref{fig:linepack}. In our setting, the gas price is approximately proportional to the gas demand. From Fig. \ref{fig:linepack}, it can be observed that the line pack plays the role of shifting the gas demand. During gas valley load periods, say period 2 to 5, extra gas are bought from the gas retailers and stored in the line pack, and during gas peak load periods, say period 9 to 11, 17 to 21, the previously stored gas are extracted from the line pack and less gas are purchased from the retailers. The percentages of stored and extracted gas in the actual gas demand are also shown in Fig. \ref{fig:linepack}. From Fig. \ref{fig:linepack}, the stored/extracted gas are comparable with the gas demand in most periods, and even larger than the gas demand, say period 3 and 4. In this way, the operation cost and operational feasibility of gas system as well as the fuel adequacy issue of gas-fired units in the power system can be great relieved if the line pack is utilized properly, validating the benefit of a multi-period OGPF model.

\begin{table}[ht]
\footnotesize
  \centering
  \caption{Parameters of Algorithm \ref{alg:ADMM} and \ref{alg:Improved_CCP}}\label{Tab:parameter}
  \begin{tabular}{c|c|c|c|c|c|c|c|c}
  \hline
  \multicolumn{3}{c|}{Algorithm \ref{alg:ADMM}}& \multicolumn{6}{c}{Algorithm \ref{alg:Improved_CCP}}\\
  \hline
   $\varsigma$&$d$ &$k_{max}$& $\delta$ & $\rho_0$ & $\rho_{max}$ & $\epsilon$ & $\kappa$ & $j^{max}$ \\
  \hline
   $10^{-3}$&100&100& 1 & 0.01 & 1000 & $10^{-6}$ & 2 & 100\\
  \hline
  \end{tabular}
\end{table}

\subsection{Comparison with other OGF Calculation Methods}

As mentioned before, there are mainly three other algebraic equation based methods to crack the nonlinear and non-convex OGF problem, which are the Nonlinear method \cite{Sheng_Probabilistic}, MILP approximation method \cite{Carlos_MILP} and the LP/SOCP relaxation method \cite{Hantao_LP,Conrado_MISOCP}. Their performances are compared from three aspects: solution quality, constraint violation, and computation time. For other OGF calculation methods, Algorithm \ref{alg:ADMM} is adopted for the OGPF problem. Particularly, in the MILP approximation method, we use an eight-segment piecewise linear approximation to replace the nonlinear Weymouth equation and an eight-segment linear approximation to represent the relationship between nodal pressure square and nodal pressure\cite{CorreaPosada2015}. The time limit is set as 7200s. Results are shown in Table \ref{Tab:Power6Gas7}.

\begin{table}[ht]
\footnotesize
  \centering
  \newcommand{\tabincell}[2]{\begin{tabular}{@{}#1@{}}#2\end{tabular}}
  \caption{Comparison with the state-of-art methodology for OGF subproblem: Power13Gas7}\label{Tab:Power6Gas7}
  \begin{tabular}{c|c|c|c|c|c}
  \hline
  \multicolumn{2}{c|}{\multirow{2}{*}{Methodology}} & Objective & \multicolumn{2}{c|}{Feasibility}& Time \\
  \cline{4-5}
  \multicolumn{2}{c|}{} & ($\times10^4$\$) & Coupling & OGF & (s) \\
  \hline
  \multicolumn{2}{c|}{Nonlinear} & $1.7639$  & $5.5\times 10^{-4}$& $9.1\times 10^{-7}$ & 82.16 \\
  \hline
  \multicolumn{2}{c|}{MILP} & $2.0464$ &  0.12 & 0.085 & $7200^{*}$\\
  \hline
  \multicolumn{2}{c|}{LP/SOCP} & $1.7022$ & $8.7\times10^{-5}$& 10.22& 0.086\\
  \hline
  \multirow{2}{*}{SSA}& SOCP & $1.7436$ & $7.2\times 10^{-5}$ & $3.1\times 10^{-7}$ & 0.574 \\
  \cline{2-6}
   & Zero & $1.8762$ & $1.9\times 10^{-4}$ & $7.2\times 10^{-7}$ & 1.544 \\
  \hline
  \end{tabular}
\end{table}

In Table \ref{Tab:Power6Gas7}, it can be observed that the MILP approximation method fails to converge in the given time limit, resulting in a suboptimal OGF solution and infeasibility of OGPF. In fact, the solution quality and computation time contradicts with each others. Using more segments will enhance the solution quality yet bring larger computational burden. In this case, eight-segment approximation causes computation overhead, which means the MILP method may not be a good choice. Other than the MILP method, all the other three methods successfully offer a solution. Among those, LP/SOCP method gives the lowest objective function in the shortest time, however, the solution is highly infeasible with a 1022\% MRCV, because the relaxation is not tight. For the proposed SSA method, the objective values and solution times with respect to different initial values differ a lot, which validates the significance of the proposed initial point selection method. Compared with the Nonlinear method, SSA is 160 times faster, and also gives a better solution. The advantages of SSA would be even more remarkable in large-scale instances.

\subsection{Computational Efficiency}

To demonstrate the scalability of the proposed method, it is applied to a larger test system, which consists of a modified IEEE 123-node power feeder and a modified Belgian high-calorific 20-node gas network. It has 6 gas-fired DGs, 4 non-gas DGs, 2 gas retailers, 3 compressors, 16 passive pipelines, 85 power loads and 9 gas loads. Please refer to \cite{Power6Gas7} for the topology, the demand curve as well as the system data. Hereinafter, the test system is referred to as Power123Gas20 system. Because the LP/SOCP relaxation method is not tight, and the MILP method is not scalable, we will focus on the comparison with the Nonlinear method. The parameters of Algorithm \ref{alg:ADMM} and \ref{alg:Improved_CCP} are the same with Table \ref{Tab:parameter} except for $d$ in Algorithm \ref{alg:ADMM}, which is $1000$ in this case. The initial point of SSA algorithm is the solution of SOCP relaxation method described in Section III.B. We also impose a 7200s CPU time limit for this case. The results are listed in Table \ref{Tab:Power39Gas20}.

\begin{table}[ht]
\footnotesize
  \centering
  \newcommand{\tabincell}[2]{\begin{tabular}{@{}#1@{}}#2\end{tabular}}
  \caption{Performance of the two algorithms: Power39Gas20}\label{Tab:Power39Gas20}
  \begin{tabular}{c|c|c|c|c}
  \hline
  \multirow{2}{*}{Periods}& \multicolumn{2}{c|}{Nonlinear} & \multicolumn{2}{c}{SSA}\\
  \cline{2-5}
  & Objective($\times 10^5\$$)& CPU(s) & Objective ($\times 10^5\$$) & CPU(s)\\
  \hline
  4 & $0.7881$  & 90.6 & $0.7775$ &3.12\\
  \hline
  8 & - & 7200*& $1.7237$ &6.97\\
  \hline
  12& - & 7200*& $2.7033$ &9.02\\
  \hline
  16& $3.6522 $ & 186  & $3.6163$ &26.3\\
  \hline
  20& -  & 7200*& $4.5393$ &7.87\\
  \hline
  24& -  & 7200*& $5.4996$ &39.7\\
  \hline
  \end{tabular}
\end{table}

From Table \ref{Tab:Power39Gas20}, SSA is significantly faster than the Nonlinear method. It is found that the Nonlinear method is not robust. It fails to converge when $T=8,12,20,24$, while SSA converges in all cases. Moreover, the CPU time of SSA is not monotonic with respect to the number of periods. One possible reason is that the quality of the initial values varies   due to the different load profile and constraint (\ref{gas_mass_storage}), which requires the line pack of each passive pipeline in the final period to be equal to that in the initial period to maintain the initial status of the next day.

\section{Conclusion}
This paper presents an SOCP and ADMM based distributed algorithm for the multi-period OGPF problem of interdependent PDN and GDN. Through the proposed algorithm, the Weymouth equation is convexified such that the OGF problem is converted to a sequence of SOCPs, which can be solved efficiently. The OGPF model underlies many engineering problems, such as economic dispatch, energy market, demand response, etc. It is interesting to study the operation flexibility brought by the power-to-gas technology and line pack, which can be regarded as a new kind of energy storage. It is also interesting to study the joint clearing of the bilateral power and natural gas market, where the locational marginal prices can be easily extracted from the Lagrangian dual multipliers at the OGPF solution.

\bibliographystyle{IEEEtran}
\bibliography{IEEEabrv,refs}

%




\end{document}